\documentclass[11pt]{article}
\usepackage{amsfonts,amssymb,amsmath,latexsym}
\usepackage[english]{babel}

\newcommand{\Rz}{\mathbb{R}}
\newcommand{\Cz}{\mathbb{C}}
\newcommand{\Kz}{\mathbb{K}}
\newcommand{\trace}{\mathrm{tr}}

\DeclareMathAlphabet{\mathbfit}{OT1}{cmr}{bx}{it}

\newenvironment {PROOF}{\textsc{Proof:} \small}{{\hspace*{\fill} $\square$}}
\newtheorem {LEMMA} {Lemma} [section]
\newtheorem {PROPOSITION} [LEMMA] {Proposition}
\newtheorem {THEOREM} [LEMMA] {Theorem}

\newtheorem {DEFINITION}[LEMMA] {Definition}
\newtheorem {REMARK}[LEMMA]{Remark}

\title {Quantum integrable Toda like systems  \\}

\author{
\textbf{Martin Bordemann\thanks{Martin.Bordemann@physik.uni-freiburg.de}~,
       \addtocounter{footnote}{2}
        Martin Walter\thanks{Martin.Walter@physik.uni-freiburg.de}} \\[3mm]
        Fakult\"at f\"ur Physik\\Universit\"at Freiburg \\
        Hermann-Herder-Str. 3 \\
        79104 Freiburg i.~Br., F.~R.~G \\[3mm]}

\date{FR-THEP-98/15 \\[1mm]
      October 1998\\[1mm]
      }

\begin {document}
\maketitle
\thispagestyle{empty}

\begin {abstract}
 Using deformation quantization and suitable 2 by 2 quantum $R$-matrices
 we show that a list of Toda like classical integrable systems given
 by Y.B.Suris is quantum integrable in the sense that the classical conserved
 quantities (which are already in involution with respect to the Poisson
 bracket) commute with respect to the standard star-product of Weyl type
 in flat $2n$-dimensional space.
\end {abstract}

\newpage

\section{Introduction}\label{IntroSec}

During the past decades a lot of new families of systems of Hamiltonian
mechanics have been found which are integrable in the sense of Liouville,
i.e. which allow for $n=\frac{1}{2}\dim ($Phasespace) independent,
Poisson-commuting integrals of motion on which the Hamiltonian functionally
depends: among these are the Calogero-Moser systems and the Toda chains
(see e.g. \cite{Per}, \cite{OPII}). These systems allow for a formulation
in terms of a Lax pair (which proves that the coefficients of the
characteristic polynomial and/or the trace polynomials of the matrix Lie
algebra valued function $L$ on phase space Poisson-commute with the
Hamiltonian) and the existence of a so-called classical $r$-matrix
expressing the Poisson bracket of the components of $L$ in terms of Lie
commutators guarantees that the above invariant functions Poisson-commute
(see e.g. \cite{bordemann1},\cite{FT}).

Already at the time of the discovery of the above-mentioned integrable
systems the question of their {\em quantum integrability}
had been considered, i.e. whether one can associate
to each of these classical conserved quantities a quantum operator
(by means of some ordering prescription) such that these quantum operators
commute with the Hamiltonian operator and among each other. In
\cite{CRM} (see also \cite{HUW92} for a much more explicit proof) arguments
were given that there was no ordering problem for
the corresponding integrals for the Calogero-Moser systems and that
the corresponding operators should commute (see also \cite{OPII} for a
similar type of argument for the Toda chain).

Recently Y.B.Suris gave a list of Toda-like systems defined by traces
over products of $2\times 2$-matrices which are all classically integrable
by means of two types of a constant classical $r$-matrix depending on
spectral parameters \cite{suris}. Among his systems are relativistic and
discretized versions of the original Toda lattice.

The motivation of this article (see also the second author's thesis
\cite{wal})
 was to check whether all the systems of
Suris' list are quantum-integrable. In order to control the possible
ordering prescriptions we chose to use the concept of deformation
quantization defined in \cite{BFFLS78} which has now been well-established
on every symplectic manifold. The advantage of this method to using
operators is the fact that the quantum noncommutative multiplication is
formulated directly on the space of classical observables as a deformed
pointwise multiplication which makes it easier and more natural to compare
with classical computations. On flat $\Rz^{2n}$ (and more generally on
every cotangent bundle, see \cite{BNWI}) there exist differential operator
representations of the deformed algebra corresponding to canonical
quantization with Weyl ordering prescription. Consequently, quantum
commutativity of the classical integrals in terms of star-products can be
translated into commuting operators if necessary. Another advantage of
deformation quantization is that quantum integrability can be formulated on
much more general symplectic manifolds where star-products still exist
thanks to the theorem of DeWilde-Lecomte \cite{DL} but operator
representations are a priori lacking.\\

We find that all the systems given by Suris are quantum integrable,
and our proof uses the quantum $R$-matrices of the form identity plus
a multiple of the classical $r$-matrix. In particular we re-obtain
the quantum integrability of the Toda chain. But, as it turned
out there is quantum asymmetry: the analogs of the coefficients of the
characteristic polynomial are actually commuting with respect to the
star-product whereas the analogs of the trace-polynomials are {\em not}.
One can cure that by adding quantum corrections to the trace polynomials
which can be obtained by interpreting the Waring identities (which
are polynomial formulas between the two sets of functions) in the deformed
algebra.

\section{Star products and ordering prescription of standard
and Weyl ordered products in $\Rz^{2n}$}

In this section we shall briefly recall the formulas needed for the
star-products and their operator representations in flat $\Rz^{2n}$
(see also \cite{AW}, \cite{BNWI}). We shall denote the co-ordinates
of $\Rz^{2n}$ by $(\vec{q},\vec{p})$.

Let $F:\Rz^{2n}\rightarrow \Cz$ be smooth. The standard ordering
prescription assigns to $F$ the formal differential operator series
\begin{equation}
     (\rho_S(F)\psi)(\vec{q}) :=
              \sum_{k=0}^\infty \frac{1}{k!}(\hbar/i)^k
          \sum_{i_1,\dots,i_k=1}^n
      \frac{\partial^k F}{\partial p_{i_1}\cdots \partial p_{i_k}}
                                                 (\vec{q},\vec{0})
      \frac{\partial^k\psi}{\partial q^{i_1}\cdots \partial q^{i_k}}(\vec{q})
\end{equation}
However, real-valued functions do not correspond to symmetric operators
(on the dense domain of compactly supported smooth complex functions),
but one rather has $\rho_S(F)^\dagger = \rho_S(N^2\bar{F})$ where
\begin{equation}
           N := exp(\frac{\hbar}{2i}\sum_{j=1}^n
                   \frac{\partial^2}{\partial q^j\partial p_j}).
\end{equation}
The Weyl ordering prescription
\begin{equation}
          \rho_W(F) := \rho_S(NF)
\end{equation}
has the more physical property $\rho_W(F)^\dagger = \rho_W(\bar{F})$
and exactly corresponds to a total symmetrization of position and momentum
operators in polynomial observables. The star-products of standard ordered
type, $\ast_S$ and of Weyl type, $\ast$, of two smooth complex-valued
functions $F,G$ on phase space are defined as follows:
\begin{equation}
  (F\ast_S G) (\vec{q},\vec{p}):=
       exp(\frac{\hbar}{i}\sum_{j=1}^n
              \frac{\partial^2}{\partial q'{}^j\partial p_i} )
       F(\vec{q},\vec{p})G(\vec{q}\,',\vec{p}\,')
              \bigg|_{\vec{q}=\vec{q}',\vec{p}=\vec{p}'}
\end{equation}
\begin{equation}
  (F\ast G) (\vec{q},\vec{p}):=
       exp(\frac{i\hbar}{2}\left( \sum_{j=1}^n
          \frac{\partial^2}{\partial q^j\partial p'_j}
              - \frac{\partial^2}{\partial {q'}{}^j\partial p_j} \right))
       F(\vec{q},\vec{p})G(\vec{q}',\vec{p}')
              \bigg|_{\vec{q}=\vec{q}',\vec{p}=\vec{p}'},
\end{equation}
they satisfy the representation identities
\begin{equation}
   \rho_S(F\ast_S G) = \rho_S(F)\rho_S(G)~~,~~
   \rho_W(F\ast G) = \rho_W(F)\rho_W(G)~~,
\end{equation}
and they are related by $N$ as follows:
\begin{equation}
     F\ast G = N^{-1}((NF)\ast_S (NG))~~.
\end{equation}
Clearly, the two star-products are associative and have the correct classical
limit, i.e. one gets pointwise multiplication at the order $\hbar^0$
and $i$ times the Poisson bracket taking the commutator at the order
$\hbar$.
For practical purposes it is often easier to compute the star-product of
standard-ordered type and to use $N$ to switch to the Weyl type
multiplication.

We conclude this section with the general definition of a quantum integrable
system:

\begin{DEFINITION}[Quantum integrable system]\label{def quant}
A classical completely integrable system
with Hamilton function $H$ on a $2n$-dimensional symplectic manifold is said
to be quantum integrable, if there exists a star product $\ast$ and
$n$ formal power series $F_i\in C^\infty(M)[[\hbar]]$ which
conincide with the classical conserved quantities $f_i$
at order zero in $\hbar$, such that
\begin{enumerate}
\item $F_i\ast H - H\ast F_i= 0$
\item $F_i\ast F_j - F_j \ast F_i= 0$
\end{enumerate}
for all $1\le i,j \le n$.
\end{DEFINITION}
\textit{Remark:}
i) Writing down the formal power serie $(1\le i\le n)$,
\[
F_i = \sum_{k=0}\hbar^k F_i^{(k)},\qquad F_i^{(0)} \equiv f_i
\]
higher terms in $\hbar$ of $F_i$ can be regarded as
quantum corrections of the classical conserved quantities $f_i$.
ii) Here, we do not consider deformations of the Hamiltonian.

\section{Quantum $R$-matrices for certain classical $r$-matrices}

The field $\Kz$ is either equal to $\Rz$ or to $\Cz$.
Let $L(n,\Kz)$ denote the space of $\Kz$-valued $(n\times n)$-matrices
with standard basis $E_{ij}$, where $1\le i,j\le n$. To express tensor
products with spectral parameters properly let
$\Kz(\lambda_1,\lambda_2,\ldots,\lambda_k)$
denote the field of rational functions over $\Kz$ in $k$ parameters
$\lambda_1,\lambda_2,\ldots,\lambda_k$. $\rho$ will be an abbreviation for
$i\hbar/2$. Furthermore, recall the standard tensor notation in the context
of $R$-matrices: for an element
$R=\sum_{i}s_i\otimes t_i \otimes \phi_i(\lambda,\mu)$ in
$L(n,\Kz)\otimes L(n,\Kz)\otimes \Kz(\lambda,\mu)$ and
for a positive integer $N\geq 2$ and integers $a,b$ s.t. $1\leq a< b \leq N$
we shall write $R_{ab}$ for
$\sum_{i}\mathbf{1}\otimes\cdots\otimes\mathbf{1}\otimes s_i\otimes
\mathbf{1}\otimes
\cdots\otimes\mathbf{1}\otimes t_i\otimes\mathbf{1}\otimes
\cdots\otimes\mathbf{1}\otimes\phi_i(\lambda_a,\lambda_b)$ (where $s_i$ is
at the $a$th tensor factor and $t_i$ is at the $b$th tensor factor) regarded
as an element in
$L(n,\Kz)^{\otimes N}\otimes\Kz(\lambda_1,\cdots,\lambda_N)$. Note that this
last space is an associative algebra with respect to tensor factor wise
multiplication. The symbol
$R_{ba}$ is equal
to the above expression with $s_i$ and $t_i$ and $\lambda_a$ and $\lambda_b$
exchanged. We shall frequently use the standard isomorphism
\begin{eqnarray}\label{iso tens prod}
L(n,\Kz) \otimes L(n,\Kz) & \to & L(n^2,\Kz)~, \nonumber \\
a\otimes b & \mapsto &
\begin{pmatrix}
a_{11}b & \dots & a_{1n}b \\
\vdots & \ddots & \vdots\\
a_{n1}b & \dots & a_{nn}b \\
\end{pmatrix}
\end{eqnarray}
to conveniently express tensor products of matrices.

We consider two particular classical $r$-matrices with spectral parameter,
i.e. elements $r$, $\tilde{r}$ in
$L(2,\Kz)\otimes L(2,\Kz)\otimes \Kz(\lambda,\mu)$ which
are antisymmetric ($r_{12}=-r_{21}$) and obey the classical Yang-Baxter
equation,
\begin{equation}\label{CYBE}
   [r_{12},r_{13}] + [r_{12},r_{23}] + [r_{13},r_{23}] = 0~~,
\end{equation}
which have been used in a preprint by Y.B.Suris \cite{suris}.
The first one is generated by the so-called
Casimir element $C:=\sum_{i,j=1}^2 E_{ij}\otimes E_{ji}$
\begin{equation}
    r := \frac{C}{\lambda-\mu}~~,
\end{equation}
the second one is
\begin{equation}
\tilde{r} = \left(
\begin{array}{cccc}
\frac{1}{2}\frac{\lambda^2+\mu^2}{\lambda^2-\mu^2} & 0 & 0 & 0\\
0 & -\frac{1}{2} & \frac{\lambda\mu}{\lambda^2-\mu^2}  & 0\\
0 &  \frac{\lambda\mu}{\lambda^2-\mu^2} & \frac{1}{2} & 0 \\
0 & 0 & 0 & \frac{1}{2}\frac{\lambda^2+\mu^2}{\lambda^2-\mu^2}
\end{array}
\right).
\end{equation}
where the above isomorphism (\ref{iso tens prod}) is used.

Define the following two quantum $R$-matrices with spectral paramter:
\begin{eqnarray}
R & = & \mathbf{1}\otimes\mathbf{1} + f(\rho)r~,\\
\tilde{R} & = & \mathbf{1}\otimes\mathbf{1} + f(\rho)\tilde{r}~,
\end{eqnarray}
where $f(\rho)$ is a smooth function and/or a formal power series.

\begin{LEMMA}
$R$ and $\tilde{R}$ fulfill the spectral
quantum Yang-Baxter equation, and are
unitary up to a factor, i.e.
\begin{eqnarray}
R_{12}R_{13}R_{23} & =    & R_{23}R_{13}R_{12} ~, \label{qybg}\\
            R_{21} & \sim & R_{12}^{-1},
\end{eqnarray}
and corresponding relations for $\tilde{R}$.
\end{LEMMA}
\begin{PROOF}
The proof is a straight forward computation using (\ref{CYBE}). It only
remains to check the vanishing of the terms of third order in $f(\rho)$.
The properties for $R$ are well-known, see e.g. \cite{CP}.
\end{PROOF}

We shall call such $R$-matrices having the properties of the
previous lemma quantum $R$-matrices in short.

\section{Quantum integrable Toda like $n$-particle systems
defined by \\($2\times 2$)-matrices}

Let $\mathcal{A}_n$ denote the associative algebra
$L(2,\Kz)\otimes C^\infty(\Rz^{2n})\otimes \Kz(\lambda)$ with tensor factor
wise multiplication and pointwise multiplication in $C^\infty(\Rz^{2n})$.
Let $\trace: \mathcal{A}_n \rightarrow
C^\infty(\Rz^{2n})\otimes \Kz(\lambda)$ denote the standard
extension of the matrix trace in the first tensor factor. For any element
$U\in \mathcal{A}_1$ let $U^k\in \mathcal{A}_n$, $1\leq k\leq n$ denote
the embedding of $\mathcal{A}_1$ into $\mathcal{A}_n$ by pulling back
the matrix elements by means of the projection $p_k:$
$\Rz^{2n}\rightarrow \Rz^2:$ $(q,p)\mapsto (q^k,p_k)$.

Given such $U\in \mathcal{A}_2$, we can built the following functions
on the phase space $\Rz^{2n}$ depending on the parameter $\lambda$:
\begin{DEFINITION}
\begin{eqnarray}
\chi_n(\lambda)         & := & \trace (U^n(\lambda)\dots U^1(\lambda)),\\
\tilde{\chi}_n(\lambda) & := & \trace (E_{11}U^n(\lambda)\dots U^1(\lambda)),
\end{eqnarray}
\end{DEFINITION}

This Definition is motivated by the following important example, the
well-known nonperiodic Toda chain
whose Hamiltonian function and Lax matrix are given by
\begin{eqnarray}
H & = & \frac{1}{2}\sum_{i=1}^n p_i^2 + \sum_{i=1}^{n-1}e^{q_i-q_{i+1}},
                                            \label{HamnonperToda}\\
L & = & \sum_{i=1}^n p_i\,E_{ii}
          + \sum_{i=1}^{n-1} e^{\frac{1}{2}(q_i-q_{i+1})}
                              \,(E_{i,i+1}+E_{i+1,i})\label{lax matrix}.
\end{eqnarray}
It is easy to see that for the choice
$U(q,p)(\lambda):=\begin{pmatrix} -\lambda+p& -e^{-q}\\ e^q & 0\end{pmatrix}$,
the above function $\tilde{\chi}_n(\lambda)$ coincides with the
characteristic polynomial $\det(-\lambda\mathbf1-L(\vec{q},-\vec{p}))$.
This particular $U$ is taken from Suris' paper \cite{suris} where our $U$
is the transposed matrix of Suris' $L$ and $\lambda$ is changed to
$-\lambda$.
The function $\chi_n$ is equal to
the characteristic polynomial of the periodic Toda chain.
It is well-known that the $n$ nonconstant coefficients
of these characteristic polynomials are functionally
independent and in involution and that the above Hamiltonian function
(\ref{HamnonperToda}) can be obtained from the coefficients of $\lambda^1$
and $\lambda^2$.

To investigate the quantum case let us define some multiplications:
\begin{DEFINITION}
 Define the star product on $\mathcal{A}_n$ as follows
 \begin{enumerate}
  \item
   \begin{eqnarray*}
     \ast:\mathcal{A}_n\times \mathcal{A}_n
                         & \to &   \mathcal{A}_n,\\
                (M \ast N)_{ik}
                         & :=  & \sum_j M_{ij}\ast N_{jk}
   \end{eqnarray*}
     where $M,N\in \mathcal{A}_n$ and $\ast$ denotes the
       standard star-product of Weyl type defined in Section 2.
  \item For $a,b\in L(2,\Kz)\otimes L(2,\Kz)$ and
         $f,g\in  C^\infty(\Rz^{2n})$ let
            $(a\otimes f)(b\otimes g):=  ab\otimes fg$.
 \end{enumerate}
\end{DEFINITION}
This definition immediately yields explicit formulae for the star product of
two characteristic polynomials.
\begin{PROPOSITION}\mbox{}

 \begin{enumerate}
  \item $\chi_n(\lambda)\ast\chi_n(\mu)
             = \trace\big(
                  (U^n_1(\lambda)\ast U^n_2(\mu)) \dots
                        (U^1_1(\lambda)\ast U^1_2(\mu)) \big)$,
  \item $\tilde{\chi}_n(\lambda)\ast\tilde{\chi}_n(\mu)
             = \trace\big( (E_{11}\otimes E_{11})
                  (U^n_1(\lambda)\ast U^n_2(\mu)) \dots
                        (U^1_1(\lambda)\ast U^1_2(\mu)) \big)$,
 \end{enumerate}
where $U_1=U\otimes \mathbf{1}$ and $U_2=\mathbf{1}\otimes U$.
\end{PROPOSITION}
\begin{PROOF}
This is shown by direct calculation. The $\lambda$- and $\mu$-dependent
matrices can be reordered without getting additional terms because the order of
terms belonging to different pairs of phase space variables is unchanged.
\end{PROOF}

The main idea to prove quantum commutativity is of course borrowed from
the theory of statistical models (see e.g. \cite{CP}) and consists in
showing that the commutation of $U(\lambda)(q,p)$
and $U(\mu)(q,p)$ as functions of the same phase space variables
can be written as a conjugation with special quantum $R$-matrices. More
precisely:
\begin{THEOREM}
 Let $U$ and $R$ be of the form as shown in the succeding table.
 Each pair of $U$ and $R$ fullfils the relation
\begin{equation}\label{rllrel}
R U_1(\lambda)\ast U_2(\mu) = U_2(\mu)\ast U_1(\lambda) R.
\end{equation}
\end{THEOREM}
\begin{center}  
\begin{tabular}{|c|c|}
\hline
Lax matrix $U(\lambda)$ & quantum $R$-matrix \\
\hline\hline
\begin{tabular}{c}
$\begin{pmatrix} -\lambda+p & -e^{-q}\\ e^{q} & 0
\end{pmatrix}$\\
$\begin{pmatrix} -\lambda +p & -e^{-q} \\  pe^q & 1  \end{pmatrix}$
\end{tabular}
&
$\mathbf{1}\otimes\mathbf{1} + 2\rho r$\\
\hline
\begin{tabular}{c}
$\begin{pmatrix} \frac{1}{\lambda} -\lambda e^p & -e^{-q} \\
e^q & -\lambda \end{pmatrix}$ \\
$\begin{pmatrix} \frac{1}{\lambda} -\lambda e^p & -g^2 e^{-q+p} \\
e^q & 0  \end{pmatrix}$\\
$\begin{pmatrix} \frac{1}{\lambda} -\lambda e^p & -g^2(e^p+\delta)e^{-q} \\
e^q & -\lambda\delta g^2 \end{pmatrix}$
\end{tabular}
&
$\mathbf{1}\otimes\mathbf{1} - 2\tanh(\rho)\tilde{r}$\\
\hline
\begin{tabular}{c}
$\begin{pmatrix}\frac{1}{\lambda} -\lambda e^{\epsilon p} & -\epsilon e^{-q} \\
\epsilon e^q & -\epsilon^2\lambda \end{pmatrix}$\\
$\begin{pmatrix}\frac{1}{\lambda} -\lambda e^{\epsilon p} & -(e^{\epsilon p}-1)e^{-q} \\
\epsilon e^q & +\lambda\epsilon \end{pmatrix}$
\end{tabular}
&
$\mathbf{1}\otimes\mathbf{1} - 2\tanh(\epsilon\rho)\tilde{r}$\\
\hline
\end{tabular}
\end{center}
\begin{PROOF}
 The assertion can directly be proved by a rather long, but straight forward
 computation.
\end{PROOF}

The presented Lax matrices $U$ originally occurred in an
`almost complete list' of Toda related classical integrable systems
formulated by $(2\times 2)$-matrices by Suris in 1997 (cf. \cite{suris}).
Note that the Lax matrices in this table differ from the ones given by Suris
by matrix transposition and the transformation
$\lambda\mapsto -\lambda$.

\begin{REMARK}
Defining the monodromy matrix $T(\lambda):=U^n(\lambda)\ast\dots
\ast U^1(\lambda)$ equation \eqref{rllrel} leads to
\begin{equation}\label{RTT}
RT_1(\lambda)\ast T_2(\mu) =  T_2(\mu)\ast T_1(\lambda)R,
\end{equation}
which is the well-known RTT-relation.
\end{REMARK}

As a corollary we obtain the main
\begin{THEOREM} We have the following quantum commutation relations:
\[
\chi_n(\lambda)\ast\chi_n(\mu) = \chi_n(\mu)\ast\chi_n(\lambda),\qquad
\tilde{\chi}_n(\lambda)\ast\tilde{\chi}_n(\mu)
      = \tilde{\chi}_n(\mu)\ast\tilde{\chi}_n(\lambda)
\]
Hence every star polynomial of coefficients of  $\chi_n(\lambda)$ or
$\tilde{\chi}_n(\lambda)$ defines a quantum integrable system. In particular,
the nonperiodic and periodic Toda chains are quantum integrable.
\end{THEOREM}

In general we can consider two sets of classical conserved
quantities: Given a Lax matrix, one takes the trace polynomials
$\{I_1,\dots,I_n\}$, where $I_k =(1/k) \trace
L^k$, or the coefficients of the characteristic polynomial
$\{J_1,\dots,J_n\}$, where we use the notation
\begin{equation}
\chi_n(\lambda) = \det(\lambda\mathbf{1}-L) = \lambda^n + \sum_{k=1}^n
J_k\lambda^{n-k}.
\end{equation}
\begin{LEMMA}
Let $r=(r_1,\dots,r_n)$ be a multiindex $(r!:=r_1!\dots r_n!$,
$|r|:=r_1+\dots+ r_n)$ and $\alpha=(1,2,\dots, n)$. The relation
between classical trace polynomials and coefficients of the characteristic
polynomial is expressed by the Waring's formulae (see e.g. \cite{waring}):
\begin{eqnarray}
J_k & = & \sum_{\alpha r=k}\frac{(-1)^{|r|}}{r!}I_1^{r_1}\dots
I_k^{r_k}, \label{j durch i}\\
I_k & = & \sum_{\alpha r=k}\frac{(-1)^{|r|}}{r!}\frac{|r|!}{|r|}J_1^{r_1}\dots J_k^{r_k}
\label{i durch j}.
\end{eqnarray}
\end{LEMMA}
Classical, the involutivity of one set follows from the
involutivity of the other set. Quantum mechanical it does not,
but, since we have shown that the ${J_k}$ $(1\le k\le n)$
are in involution, we know that the star commutator of the
corresponding star polynomials of \eqref{i durch j}, say
$\hat{I}_k$, which one gets by replacing the usual multiplication
by star products, vanishes. In the lowest order in $\hbar$ $\hat{I}_k$
coincides with $I_k$. Therefore we can calculate quantum
corrections for the classical trace polynomials of the Toda chain
with Lax matrix \eqref{lax matrix}: For $k\le 3$ $\hat{I}_k$ and $I_k$
are the same, for greater $k$ one gets:
\begin{eqnarray}
\hat{I}_4  & = &  I_4 + \rho^2 \sum_{i=1}^{n-1} e^{q_i-q_{i+1}},\\
\hat{I}_5  & = &  I_5 + 2\rho^2\sum_{i=1}^{n-1} (p_i+p_{i+1}) e^{q_i-q_{i+1}},\\
\hat{I}_6  & = &  I_6 + \rho^2\sum_{i=1}^{n-1}\left(\frac{8}{3}e^{2(q_i-q_{i+1})}
                   + \frac{10}{3}(p_i^2 + p_ip_{i+1} + p_{i+1}^2)e^{q_i-q_{i+1}}
                   \right)\nonumber \\
& & \quad          + \rho^2 \sum_{i=1}^{n-2}e^{q_i-q_{i+2}} +\rho^4
                   \sum_{i=1}^{n-1}e^{q_i-q_{i+1}}.
\end{eqnarray}

\section*{Acknowledgment}
We would like to thank P. Kulish, J. M. Maillet and
S. Waldmann for various discussions.

\end {document}